\documentclass[12pt]{article}
\usepackage{mathrsfs}
\usepackage{epic,eepic,epsf,epsfig}
\usepackage{amsfonts,srcltx,mathrsfs}
\usepackage{multirow}
\usepackage{amsmath}
\usepackage{amssymb}
\usepackage{amsbsy}
\usepackage{graphicx}
\usepackage{amsfonts}
\usepackage{color}
\usepackage{setspace}

\newtheorem{lem}{Lemma}[section]%
\newtheorem{theorem}[lem]{Theorem}%
\newtheorem{cor}[lem]{Corollary}%

\def\nd{\mathrel{\bigm|\kern-.7em/}}

\def\f{\noindent}

\def\P\GammaL{\hbox{\rm P\GammaL}}

\begin{document}
\title{Walks, infinite series and spectral radius of graphs}

\footnotetext{E-mails: zhangwq@pku.edu.cn}

\author{Wenqian Zhang\\
{\small School of Mathematics and Statistics, Shandong University of Technology}\\
{\small Zibo, Shandong 255000, P.R. China}}
\date{}
\maketitle

\begin{abstract}
For a graph $G$, the spectral radius $\rho(G)$ of $G$ is the largest eigenvalue of its adjacency matrix. In this paper, we seek the relationship between $\rho(G)$ and the walks of the subgraphs of $G$. Especially, if $G$ contains a complete multi-partite graph as a spanning subgraph, we give a formula for $\rho(G)$ by using an infinite series on walks of the subgraphs of $G$. These results are useful for the current popular spectral extremal problem.

\bigskip

\f {\bf Keywords:} spectral radius; spectral extrema; walk; infinite series.

\end{abstract}

\baselineskip 17 pt

\section{Introduction}

All graphs considered in this paper are finite, undirected and without parallel edges or loops.
For a graph $G$, let $\overline{G}$ be its complement. The vertex set and edge set of $G$ are denoted by $V(G)$ and $E(G)$, respectively. For a vertex $u$, let $d_{G}(u)$ be its degree. For two vertices $u$ and $v$, we say $u\sim v$ if they are adjacent in $G$.  For a certain integer $n$, let $K_{n}$ and $K_{1,n-1}$ be the complete graph and the star of order $n$, respectively.  For $\ell\geq2$ vertex-disjoint graphs $G_{1},G_{2},...,G_{\ell}$, let $\cup_{1\leq i\leq \ell}G_{i}$ be the disjoint union of them. Let $G_{1}\vee G_{2}$ be the join of $G_{1}$ and $G_{2}$. For any terminology used but not defined here, one may refer to \cite{CRS}.

 Let $G$ be a graph with vertices $v_{1},v_{2},...,v_{n}$. The {\em adjacency matrix} of $G$ is an $n\times n$ matrix $(a_{ij})$, where $a_{ij}=1$ if  $v_{i}\sim v_{j}$, and $a_{ij}=0$ otherwise. The {\em spectral radius} $\rho(G)$ of $G$ is the largest eigenvalue of its adjacency matrix. By Perron--Frobenius theorem,  $\rho(G)$ has a non-negative eigenvector.  Any non-negative eigenvector corresponding to $\rho(G)$ is called a {\em Perron vector} of $G$. If $G$ is connected, any Perron vector of $G$ has positive entries.

For a family $\mathcal{H}$ of graphs, a graph $G$ is call $\mathcal{H}$-free if $G$ does not contain any one in $\mathcal{H}$ as a subgraph. Let ${\rm EX}(n,\mathcal{H})$ be the set of $\mathcal{H}$-free graphs of order $n$ with the maximum number of edges, and let ${\rm SPEX}(n,\mathcal{H})$ be the set of $\mathcal{H}$-free graphs of order $n$ with the maximum spectral radius. For $r\geq2,n\geq r$, let $T_{n,r}$ be the Tur\'{a}n graph of order $n$ with $r$ parts. Tur\'{a}n's famous theorem states that ${\rm EX}(n,\left\{K_{r+1}\right\})=\left\{T_{n,r}\right\}$.
Nikiforov \cite{N2} proves that ${\rm SPEX}(n,\left\{K_{r+1}\right\})=\left\{T_{n,r}\right\}$.
 In 2010, Nikiforov \cite{N1} formally presented a spectral version of Tur\'{a}n-type
problem: what is the maximum spectral radius of an $\mathcal{H}$-free graph of order $n$?
 In recent years, this problem has become very popular (for example, see \cite{BDT,CFTZ,CDT1,CDT2,CDT3,DKLNTW,FLSZ,FTZ,LL,LP1,LP2,LZZ,NWK,WKX,ZL}). In particular, Wang, Kang and Xue \cite{WKX} proved the following nice result, which is conjectured by Cioab\u{a}, Desai and Tait \cite{CDT2}.

\begin{theorem}{\rm (Wang, Kang and Xue \cite{WKX})} \label{Wang Kang Xue}
For $r\geq2$ and sufficiently large $n$, let $\mathcal{H}$ be a finite family of graphs such that all
graphs in ${\rm EX}(n,\mathcal{H})$ are obtained from $T_{n,r}$ by embedding $k$ edges into its $r$ parts, where $k$ is a bounded integer. Then ${\rm SPEX}(n,\mathcal{H})\subseteq {\rm EX}(n,\mathcal{H})$.
\end{theorem}

In general, ${\rm SPEX}(n,\mathcal{H})$ is a proper subset of ${\rm EX}(n,\mathcal{H})$ in Theorem \ref{Wang Kang Xue}. To precisely characterize the graphs in ${\rm SPEX}(n,\mathcal{H})$, we need to develop new techniques.

Let $G$ be a graph. For an integer $\ell\geq1$, a {\em walk} of length $\ell$ in $G$ is an ordered sequence of vertices $v_{0},v_{1},...,v_{\ell}$, such that $v_{i}\sim v_{i+1}$ for any $0\leq i\leq \ell-1$. The walk is called an $\ell$-walk, and the vertices $v_{0},v_{\ell}$ are called the starting and end vertices of the walk respectively. For a vertex $u$ of $G$, let $w^{\ell}_{G}(u)$ be the number of $\ell$-walks starting at $u$ in $G$. Let $W^{\ell}(G)=\sum_{v\in V(G)}w^{\ell}_{G}(u)$.
Clearly, $w^{1}_{G}(u)=d_{G}(u)$ and $$w^{\ell+1}_{G}(u)=\sum_{v\in V(G),v\sim u}w^{\ell}_{G}(v)$$
 for any $\ell\geq1$.
Thus $W^{1}(G)=2|E(G)|$,
and $$W^{2}(G)=\sum_{u\in V(G)}\sum_{v\in V(G),v\sim u}d_{G}(v)=\sum_{v\in V(G)}d^{2}_{G}(v).$$

\begin{theorem}\label{multi-set}
For an integer $r\geq2$, let $K_{n_{1},n_{2},...,n_{r}}$ be the complete $r$-partite graph of order $n$ with parts $V_{1},V_{2},...,V_{r}$, where $n_{i}=|V_{i}|\geq1$ for $1\leq i\leq r$ and $\sum_{1\leq i\leq r}n_{i}=n$. For each $1\leq i\leq r$, let $H_{i}$ be a graph with vertex set $V(H_{i})\subseteq V_{i}$. Let $G$ be the graph obtained from $K_{n_{1},n_{2},...,n_{r}}$ by embedding the edges of $H_{i}$ into $V_{i}$ for all $1\leq i\leq r$. Set $\rho=\rho(G)$. Then $$\sum_{1\leq s\leq r}\frac{1}{1+\frac{n_{s}}{\rho}+\sum^{\infty}_{i=1}\frac{W^{i}(H_{s})}{\rho^{i+1}}}=r-1.$$
\end{theorem}

Let $G$ be a graph and $B\subseteq V(G)$. For a vertex $u$ of $G$, let $w^{\ell}_{G,B}(u)$ be the number of $\ell$-walks of $G$ starting at $u$ and end at some vertex in $B$. Let $W^{\ell}_{B}(G)=\sum_{v\in B}w^{\ell}_{G,B}(v)$.

\begin{theorem}\label{fish-type}
Let $G$ be a connected graph of order $n$, and let $S\subseteq V(G)$ with $1\leq|S|<n$ such that $G[S]$ is $n_{1}$-regular. Assume that $V(G)-S=B\cup T$, where each vertex in $B$ is adjacent to all the vertices in $S$, and each vertex in $T$ is adjacent to none of the vertices in $S$. Set $H=G-S$ and $\rho=\rho(G)$. Then 
$$\frac{|S|}{\rho(\rho-n_{1})}\left(|B|+\sum_{i\geq1}\frac{W^{i}_{B}(H)}{\rho^{i}}\right)=1.$$
\end{theorem}

Let $G_{1}$ and $G_{2}$ be two graphs (probably of different orders). We say that $G_{1}$ is strictly {\em walk-preferable} to $G_{2}$ (denoted by $G_{1}\succ G_{2}$), if there is an integer $\ell\geq1$ such that $W^{\ell}(G_{1})>W^{\ell}(G_{2})$ and $W^{i}(G_{1})=W^{i}(G_{2})$ for any $1\leq i\leq \ell-1$;  $G_{1}$ is {\em walk-equivalent} to $G_{2}$ (denoted by $G_{1}\equiv G_{2}$), if $W^{i}(G_{1})=W^{i}(G_{2})$ for any $i\geq1$. We say $G_{1}$ is {\em walk-preferable} to $G_{2}$ (denoted by $G_{1}\succeq G_{2}$), if $G_{1}\succ G_{2}$ or $G_{1}\equiv G_{2}$. We also write $G_{2}\prec G_{1}$ or $G_{2}\preceq G_{1}$, if $G_{1}\succ G_{2}$ or $G_{1}\succeq G_{2}$. Let $\mathcal{G}$ be a set of graphs. Let ${\rm SPEX}(\mathcal{G})$ denote the set of graphs with the maximum spectral radius in $\mathcal{G}$.  Define ${\rm EX}^{1}(\mathcal{G})=\left\{G\in\mathcal{G}~|~W^{1}(G)\geq W^{1}(G')~ {\rm for~any}~G'\in\mathcal{G}\right\}$, and for any $\ell\geq2$ define ${\rm EX}^{\ell}(\mathcal{G})=\left\{G\in{\rm EX}^{\ell-1}(\mathcal{G})~|~W^{\ell}(G)\geq W^{\ell}(G')~ {\rm for~any}~G'\in{\rm EX}^{\ell-1}(\mathcal{G})\right\}$. Clearly, ${\rm EX}^{1}(\mathcal{G})$ is precisely the set of graphs in $\mathcal{G}$ with the maximum number of edges, and ${\rm EX}^{i+1}(\mathcal{G})\subseteq {\rm EX}^{i}(\mathcal{G})$ for any $i\geq1$. Let ${\rm EX}^{\infty}(\mathcal{G})=\cap_{1\leq i\leq\infty}{\rm EX}^{i}(\mathcal{G})$. The graphs in ${\rm EX}^{\infty}(\mathcal{G})$ are called the most walk-preferable graphs in $\mathcal{G}$.

\medskip

\begin{theorem}\label{one-set}
Let $G$ be a connected graph of order $n$, and let $S$ be a subset of $V(G)$ with $1\leq|S|<n$. Assume that $T$ is a set of some isolated vertices of $G-S$, such that each vertex in $T$ is adjacent to each vertex in $S$ in $G$. Let $H_{1}$ and $H_{2}$ be two graphs with vertex set $T$. For any $1\leq i\leq2$, let $G_{i}$ be the graph obtained from $G$ by embedding the edges of $H_{i}$ into $T$. When $\rho(G)$ is sufficiently large (compared with $|T|$), we have the following conclusions.\\
$(i)$ If $H_{1}\equiv H_{2}$, then $\rho(G_{1})=\rho(G_{2})$.\\
$(ii)$ If $H_{1}\succ H_{2}$, then $\rho(G_{1})>\rho(G_{2})$.\\
$(iii)$ If $H_{1}\prec H_{2}$, then $\rho(G_{1})<\rho(G_{2})$.
\end{theorem}

The rest of the paper is organized as follows. In Section 2, we give the proofs of Theorems \ref{multi-set}, \ref{fish-type} and \ref{one-set}. In Section 3, some applications will be included.

\section{Proofs of Theorems \ref{multi-set}, \ref{fish-type} and \ref{one-set}}

The following Lemma \ref{subgraph} is taken  from Theorem 8.1.1 of \cite{CRS}.

\begin{lem}{\rm (\cite{CRS})}\label{subgraph}
If $H$ is a subgraph of a connected graph $G$, then $\rho(H)\leq\rho(G)$, with equality if and only if $H=G$.
\end{lem}

\begin{lem}\label{exist}
Let $G$ be a graph of order $n$.  Then $\sum^{\infty}_{k=1}\frac{W^{k}(G)}{x^{k}}$ exists for any $x>\rho(G)$.
\end{lem}

\f{\bf Proof:}  Let $A$ be the adjacency matrix of $G$, and let $\lambda_{1}\geq\lambda_{2}\geq\cdots \geq\lambda_{n}$ be the eigenvalues of $G$. Assume that $\mathbf{a}_{i}$ is a unit eigenvector corresponding to $\lambda_{i}$ for any $1\leq i\leq n$, and $\mathbf{a}_{1},\mathbf{a}_{2},...,\mathbf{a}_{n}$ are mutually orthogonal. Let $\mathbf{j}$ be the all-one vector of dimension $n$. By linear algebra, we have $A^{k}=\sum_{1\leq i\leq n}\lambda^{k}_{i}\mathbf{a}_{i}\mathbf{a}'_{i}$ for any $k\geq1$, where $\mathbf{a}'_{i}$ is the transpose of $\mathbf{a}_{i}$.
It follows that
 $$W^{k}(G)=\mathbf{j}'A^{k}\mathbf{j}=\sum_{1\leq i\leq n}\lambda^{k}_{i}\mathbf{j}'\mathbf{a}_{i}\mathbf{a}'_{i}\mathbf{j}=\sum_{1\leq i\leq n}\lambda^{k}_{i}(\mathbf{a}'_{i}\mathbf{j})^{2}.$$
  Set $a_{i}=\mathbf{a}'_{i}\mathbf{j}$ for any $1\leq i\leq n$. Then $W^{k}(G)=\sum_{1\leq i\leq n}\lambda^{k}_{i}a^{2}_{i}$ for any $k\geq1$.
Since $x>\rho(G)$, we see that $|\frac{\lambda_{i}}{x}|<1$  for any $1\leq i\leq n$. It follows that
$$\sum^{\infty}_{k=1}\frac{W^{k}(G)}{x^{k}}=\sum^{\infty}_{k=1}\frac{\sum_{1\leq i\leq n}\lambda^{k}_{i}a^{2}_{i}}{x^{k}}=\sum_{1\leq i\leq n}a^{2}_{i}\sum^{\infty}_{k=1}(\frac{\lambda_{i}}{x})^{k}$$
exists. This completes the proof.  \hfill$\Box$

\begin{lem}\label{entry}
Let $G$ be a connected graph of order $n$, and let $S$ be a subset of $V(G)$ with $1\leq|S|<n$. Assume that $H$ is a component of $G-S$, such that each vertex of $H$ is adjacent to each vertex in $S$. Let $\rho=\rho(G)$ and let $\mathbf{x}=(x_{v})_{v\in V(G)}$ be the Perron vector of $G$ such that $\sum_{v\in S}x_{v}=\rho$. Then for any $u\in V(H)$ we have $$x_{u}=1+\sum^{\infty}_{i=1}\frac{w_{H}^{i}(u)}{\rho^{i}}.$$
\end{lem}

\f{\bf Proof:}  Denote $V(H)=T$. Recall that $\sum_{v\in S}x_{v}=\rho$. For any $u\in T$,
we have $$\rho x_{u}=(\sum_{v\in S}x_{v})+\sum_{v\sim u,v\in T} x_{v}=\rho+\sum_{v\sim u,v\in T} x_{v}.$$
 This implies that $x_{u}=1+\frac{\sum_{v\sim u,v\in T} x_{v}}{\rho}$.

\medskip

\f{\bf Claim 1.} For any vertex $u$ in $T$, we have $x_{u}\geq 1+\sum_{1\leq i\leq\ell}\frac{w_{H}^{i}(u)}{\rho^{i}}$ for any integer $\ell\geq1$.

\medskip

\f{\bf Proof of Claim 1.}
From $x_{u}=1+\frac{\sum_{v\sim u,v\in T} x_{v}}{\rho}$ we see that $x_{u}\geq1$ for any $u\in T$.
 Thus, for any $v\in T$ and $v\sim u$ we have $ x_{v}\geq1$. By $x_{u}=1+\frac{\sum_{v\sim u,v\in T} x_{v}}{\rho}$ we have
$x_{u}\geq1+\frac{w_{H}^{1}(u)}{\rho}$ for any $u\in T$. Thus Claim 1 holds for $\ell=1$.
Assume that Claim 1 holds for $\ell=j\geq1$. That is $ x_{u}\geq1+\sum_{1\leq i\leq j}\frac{w_{H}^{i}(u)}{\rho^{i}}$ for any $u\in T$.
Now we consider the case $\ell= j+1$. Let $u$ be a vertex in $T$. Note that $ x_{v}\geq1+\sum_{1\leq i\leq j}\frac{w_{H}^{i}(v)}{\rho^{i}}$ for any $v\in T$ and $v\sim u$ by induction on $\ell=j$. By $x_{u}=1+\frac{\sum_{v\sim u, v\in T} x_{v}}{\rho}$, we have
$$x_{u}\geq1+\frac{\sum_{v\sim u,v\in T} (1+\sum_{1\leq i\leq j}\frac{w_{H}^{i}(v)}{\rho^{i}})}{\rho}=1+\sum_{1\leq i\leq j+1}\frac{w_{H}^{i}(u)}{\rho^{i}}.$$
So, Claim 1 holds for $\ell=j+1$. By induction on $\ell$, we complete the proof of Claim 1. \hfill$\Box$

\medskip

Let $u^{*}$ be a vertex in $T$ such that $x_{u^{*}}=\max_{v\in T} x_{v}$.

\medskip

\f{\bf Claim 2.} For any vertex $u$ in $T$, we have $x_{u}\leq(1+\sum_{1\leq i\leq\ell}\frac{w_{H}^{i}(u)}{\rho^{i}})+(x_{u^{*}}-1)\frac{w_{H}^{\ell}(u)}{\rho^{\ell}}$ for any integer $\ell\geq1$.

\medskip

\f{\bf Proof of Claim 2.}  
By Claim 1,  we see that $x_{u^{*}}\geq1$.
For any vertex $u\in T$, from $x_{u}=1+\frac{\sum_{v\sim u,v\in T} x_{v}}{\rho}$ we obtain that
 $$x_{u}\leq1+\frac{w_{H}^{1}(u)\cdot x_{u^{*}}}{\rho}=(1+\frac{ w_{H}^{1}(u)}{\rho})+(x_{u^{*}}-1)\frac{ w_{H}^{1}(u)}{\rho}.$$
 Thus Claim 2 holds for $\ell=1$. Assume that Claim 2 holds for $\ell=j\geq1$. That is $x_{u}\leq(1+\sum_{1\leq i\leq j}\frac{w_{H}^{i}(u)}{\rho^{i}})+(x_{u^{*}}-1)\frac{w_{H}^{j}(u)}{\rho^{j}}$ for any $u\in T$.
Now we consider the case $\ell=j+1$. Let $u\in T$. Note that $ x_{v}\leq(1+\sum_{1\leq i\leq j}\frac{w_{H}^{i}(v)}{\rho^{i}})+(x_{u^{*}}-1)\frac{w_{H}^{j}(v)}{\rho^{j}}$ for any $v\in T$ and $v\sim u$ by induction on $\ell=j$. By $x_{u}=1+\frac{\sum_{v\sim u, v\in T} x_{v}}{\rho}$, we have
$$x_{u}\leq1+\frac{\sum_{v\sim u,v\in T} ((1+\sum_{1\leq i\leq j}\frac{w_{H}^{i}(v)}{\rho^{i}})+(x_{u^{*}}-1)\frac{w_{H}^{j}(v)}{\rho^{j}})}{\rho}=(1+\sum_{1\leq i\leq j+1}\frac{w_{H}^{i}(u)}{\rho^{i}})+(x_{u^{*}}-1)\frac{w_{H}^{j+1}(u)}{\rho^{j+1}}.$$
So, Claim 2 holds for $\ell=j+1$. By induction on $\ell$, we complete the proof of Claim 2. \hfill$\Box$

\medskip

By Lemma \ref{subgraph}, we have $\rho>\rho(H)$.  By Lemma \ref{exist}, $\sum^{\infty}_{i=1}\frac{W^{i}(H)}{\rho^{i}}$ exists. Let $u\in T$. Then   $\sum^{\infty}_{i=1}\frac{w_{H}^{i}(u)}{\rho^{i}}$ exists. This implies that $\frac{w_{H}^{\ell}(u)}{\rho^{\ell}}\rightarrow0$ for $\ell\rightarrow\infty$.
By Claim 1 and Claim 2 we have $$1+\sum_{1\leq i\leq\ell}\frac{w_{H}^{i}(u)}{\rho^{i}}\leq x_{u}\leq(1+\sum_{1\leq i\leq\ell}\frac{w_{H}^{i}(u)}{\rho^{i}})+(x_{u^{*}}-1)\frac{w_{H}^{\ell}(u)}{\rho^{\ell}},$$
for any integer $\ell\geq1$.
By letting $\ell\rightarrow\infty$ we obtain  that $x_{u}=1+\sum^{\infty}_{i=1}\frac{w_{H}^{i}(u)}{\rho^{i}}$.
This completes the proof. \hfill$\Box$

\medskip

\medskip

\f{\bf Proof of Theorem \ref{multi-set}.}  Let $\mathbf{x}=(x_{v})_{v\in V(G)}$ be the Perron vector of $G$ such that $\sum_{v\in V(G)}x_{v}=1$. For each $1\leq s\leq r$, let $\mathbf{y}^{s}=\frac{\rho}{1-\sum_{v\in V_{s}}x_{v}}\mathbf{x}$. Then $\mathbf{y}^{s}=(y^{s}_{v})_{v\in V(G)}$ is a Perron vector of $G$ such that $$\sum_{v\in V(G)-V_{s}}y^{s}_{v}=\frac{\rho}{1-\sum_{v\in V_{s}}x_{v}}\sum_{v\in V(G)-V_{s}}x_{v}=\rho.$$
By Lemma \ref{entry}, we have $$y^{s}_{u}=1+\sum^{\infty}_{i=1}\frac{w_{H_{s}}^{i}(u)}{\rho^{i}},$$
for any $u\in V_{s}$.
It follows that
$$x_{u}=\frac{1-\sum_{v\in V_{s}}x_{v}}{\rho}(1+\sum^{\infty}_{i=1}\frac{w_{H_{s}}^{i}(u)}{\rho^{i}}),$$
for any $u\in V_{s}$.
 Thus $$\sum_{u\in V_{s}}x_{u}=\frac{1-\sum_{v\in V_{s}}x_{v}}{\rho}\sum_{u\in V_{s}}(1+\sum^{\infty}_{i=1}\frac{w_{H_{s}}^{i}(u)}{\rho^{i}})=\frac{1-\sum_{u\in V_{s}}x_{u}}{\rho}(n_{s}+\sum^{\infty}_{i=1}\frac{W^{i}(H_{s})}{\rho^{i}}).$$
It follows that
 $$\sum_{u\in V_{s}}x_{u}=1-\frac{1}{1+\frac{n_{s}}{\rho}+\sum^{\infty}_{i=1}\frac{W^{i}(H_{s})}{\rho^{i+1}}}.$$
 Then
  $$1=\sum_{1\leq s\leq r}\sum_{u\in V_{s}}x_{u}=r-\sum_{1\leq s\leq r}\frac{1}{1+\frac{n_{s}}{\rho}+\sum^{\infty}_{i=1}\frac{W^{i}(H_{s})}{\rho^{i+1}}},$$
 implying that
 $$\sum_{1\leq s\leq r}\frac{1}{1+\frac{n_{s}}{\rho}+\sum^{\infty}_{i=1}\frac{W^{i}(H_{s})}{\rho^{i+1}}}=r-1.$$
This completes the proof. \hfill$\Box$

\medskip

\f{\bf Proof of Theorem \ref{fish-type}.} Note $V(H)=B\cup T$. Each vertex in $B$ is adjacent to all the vertices in $S$, and the vertices in $T$ has no neighbors in $S$.
Let $\mathbf{x}=(x_{v})_{v\in V(G)}$ be the Perron vector of $G$ such that $\sum_{v\in S}x_{v}=\rho$. 
Using a very similar discussion to Lemma \ref{entry}, we can show that
$$x_{u}=1+\sum^{\infty}_{i=1}\frac{w_{H,B}^{i}(u)}{\rho^{i}},$$
for any $u\in V(H)$. It follows that
$$\sum_{v\in B}x_{v}=|B|+\sum^{\infty}_{i=1}\sum_{v\in B}\frac{w_{H,B}^{i}(v)}{\rho^{i}}
=|B|+\sum^{\infty}_{i=1}\frac{W_{B}^{i}(H)}{\rho^{i}}.$$

Recall that $G[S]$ is $n_{1}$-regular. Then 
$$\rho\sum_{v\in S}x_{v}=n_{1}\sum_{v\in S}x_{v}+|S|\sum_{v\in B}x_{v}.$$
Thus, 
$$\sum_{v\in B}x_{v}=\frac{\rho-n_{1}}{|S|}\sum_{v\in S}x_{v}=\frac{\rho(\rho-n_{1})}{|S|}.$$
Consequently,
$$\frac{|S|}{\rho(\rho-n_{1})}\left(|B|+\sum_{i\geq1}\frac{W^{i}_{B}(H)}{\rho^{i}}\right)=1.$$
This completes the proof. \hfill$\Box$

\medskip

\f{\bf Proof of Theorem \ref{one-set}}.  Clearly, both $G_{1}$ and $G_{2}$ are connected. Let $\rho_{1}=\rho(G_{1})$ and $\rho_{2}=\rho(G_{2})$. Then $\rho_{1},\rho_{2}\geq\rho(G)$ by Lemma \ref{subgraph}. Let $\mathbf{x}=(x_{v})_{v\in V(G_{1})}$ be a Perron vector of $G_{1}$ such that $\sum_{v\in S}x_{v}=\rho_{1}$.  Let $\mathbf{y}=(y_{v})_{v\in V(G_{2})}$ be a Perron vector of $G_{2}$ such that $\sum_{v\in S}x_{v}=\rho_{2}$.

\medskip

\f{\bf Claim 1.} $$(\rho_{1}-\rho_{2})g(\rho_{1},\rho_{2})
=\sum^{\infty}_{\ell=0}(\frac{W^{\ell+1}(H_{1})}{\rho_{1}^{\ell}}-\frac{W^{\ell+1}(H_{2})}{\rho_{2}^{\ell}}),$$
where $$g(\rho_{1},\rho_{2})=\left(\sum_{v\in V(G)}x_{v}y_{v}\right)-\left(\sum_{u\in T}\sum^{\infty}_{\ell=1}\sum_{1\leq i\leq \ell}\frac{w^{i}_{H_{1}}(u)w^{\ell+1-i}_{H_{2}}(u)}{\rho^{i}_{1}\rho^{\ell+1-i}_{2}}\right)\geq|T|>0.$$

\medskip

\f{\bf Proof of Claim 1.}
Clearly, the only distinct edges of $G_{1}$ and $G_{2}$ are inside $T$. Let $A(G_{1})$ and $A(G_{2})$ be the adjacency  matrices of $G_{1}$ and $G_{2}$, respectively. Let $\mathbf{x}^{'}$ be the transpose of  $\mathbf{x}$. For any $u\in T$, we write $w^{0}_{H_{i}}(u)=1$ for any $1\leq i\leq2$. Thus $x_{u}=\sum^{\infty}_{i=0}\frac{w_{H_{1}}^{i}(u)}{\rho_{1}^{i}}$ and $y_{u}=\sum^{\infty}_{i=0}\frac{w_{H_{2}}^{i}(u)}{\rho_{2}^{i}}$ by Lemma \ref{entry}.
Then
\begin{equation}
\begin{aligned}
&(\rho_{1}-\rho_{2})\sum_{v\in V(G)}x_{v}y_{v}\\
&=\mathbf{x}^{'}(A(G_{1})-A(G_{2}))\mathbf{y}\\
&=\left(\sum_{uv\in E(H_{1})}(x_{u}y_{v}+x_{v}y_{u})\right)-\left(\sum_{uv\in E(H_{2})}(x_{u}y_{v}+x_{v}y_{u})\right)\\
&=\left(\sum_{u\in T}y_{u}\sum_{v\in T,vu\in E(H_{1})}x_{v}\right)-\left(\sum_{u\in T}x_{u}\sum_{v\in T,vu\in E(H_{2})}y_{v}\right)\\
&=\left(\sum_{u\in T}(\sum^{\infty}_{i=0}\frac{w_{H_{2}}^{i}(u)}{\rho_{2}^{i}})(\sum_{v\in T,vu\in E(H_{1})}\sum^{\infty}_{i=0}\frac{w_{H_{1}}^{i}(v)}{\rho_{1}^{i}})\right)
-\left(\sum_{u\in T}(\sum^{\infty}_{i=0}\frac{w_{H_{1}}^{i}(u)}{\rho_{1}^{i}})(\sum_{v\in T,vu\in E(H_{2})}\sum^{\infty}_{i=0}\frac{w_{H_{2}}^{i}(v)}{\rho_{2}^{i}})\right)\\
&=\left(\sum_{u\in T}(\sum^{\infty}_{i=0}\frac{w_{H_{2}}^{i}(u)}{\rho_{2}^{i}})(\sum^{\infty}_{i=0}\frac{w_{H_{1}}^{i+1}(u)}{\rho_{1}^{i}})\right)
-\left(\sum_{u\in T}(\sum^{\infty}_{i=0}\frac{w_{H_{1}}^{i}(u)}{\rho_{1}^{i}})(\sum^{\infty}_{i=0}\frac{w_{H_{2}}^{i+1}(u)}{\rho_{2}^{i}})\right)\\
&=\left(\sum_{u\in T}\sum^{\infty}_{\ell=0}\sum_{i,j\geq0,i+j=\ell}\frac{w_{H_{1}}^{i+1}(u)w_{H_{2}}^{j}(u)}{\rho_{1}^{i}\rho_{2}^{j}}\right)
-\left(\sum_{u\in T}\sum^{\infty}_{\ell=0}\sum_{i,j\geq0,i+j=\ell}\frac{w_{H_{1}}^{i}(u)w_{H_{2}}^{j+1}(u)}{\rho_{1}^{i}\rho_{2}^{j}}\right)\\
&=\left(\sum_{u\in T}\sum^{\infty}_{\ell=0}(\frac{w_{H_{1}}^{\ell+1}(u)w_{H_{2}}^{0}(u)}{\rho_{1}^{\ell}\rho_{2}^{0}}-
\frac{w_{H_{1}}^{0}(u)w_{H_{2}}^{\ell+1}(u)}{\rho_{1}^{0}\rho_{2}^{\ell}})\right)+\left(\sum_{u\in T}\sum^{\infty}_{\ell=1}(\rho_{1}-\rho_{2})\sum_{1\leq i\leq \ell}\frac{w^{i}_{H_{1}}(u)w^{\ell+1-i}_{H_{2}}(u)}{\rho^{i}_{1}\rho^{\ell+1-i}_{2}}\right)\\
&=\left(\sum^{\infty}_{\ell=0}(\frac{W^{\ell+1}(H_{1})}{\rho_{1}^{\ell}}-\frac{W^{\ell+1}(H_{2})}{\rho_{2}^{\ell}})\right)
+(\rho_{1}-\rho_{2})\left(\sum_{u\in T}\sum^{\infty}_{\ell=1}\sum_{1\leq i\leq \ell}\frac{w^{i}_{H_{1}}(u)w^{\ell+1-i}_{H_{2}}(u)}{\rho^{i}_{1}\rho^{\ell+1-i}_{2}}\right).
\end{aligned}\notag
\end{equation}
It follows that
$$(\rho_{1}-\rho_{2})g(\rho_{1},\rho_{2})
=\sum^{\infty}_{\ell=0}(\frac{W^{\ell+1}(H_{1})}{\rho_{1}^{\ell}}-\frac{W^{\ell+1}(H_{2})}{\rho_{2}^{\ell}}),$$
where
$$g(\rho_{1},\rho_{2})=\left(\sum_{v\in V(G)}x_{v}y_{v}\right)-\left(\sum_{u\in T}\sum^{\infty}_{\ell=1}\sum_{1\leq i\leq \ell}\frac{w^{i}_{H_{1}}(u)w^{\ell+1-i}_{H_{2}}(u)}{\rho^{i}_{1}\rho^{\ell+1-i}_{2}}\right).$$
Note that $$\sum_{v\in V(G)}x_{v}y_{v}\geq\sum_{u\in T}x_{u}y_{u}=\sum_{u\in T}(\sum^{\infty}_{i=0}\frac{w^{i}_{H_{1}}(u)}{\rho^{i}_{1}})(\sum^{\infty}_{j=0}\frac{w^{j}_{H_{2}}(u)}{\rho^{j}_{2}}).$$
Thus
\begin{equation}
\begin{aligned}
g(\rho_{1},\rho_{2})
&\geq\left(\sum_{u\in T}(\sum^{\infty}_{i=0}\frac{w^{i}_{H_{1}}(u)}{\rho^{i}_{1}})(\sum^{\infty}_{j=0}\frac{w^{j}_{H_{2}}(u)}{\rho^{j}_{2}})\right)
-\left(\sum_{u\in T}\sum^{\infty}_{\ell=1}\sum_{1\leq i\leq \ell}\frac{w^{i}_{H_{1}}(u)w^{\ell+1-i}_{H_{2}}(u)}{\rho^{i}_{1}\rho^{\ell+1-i}_{2}}\right)\\
&=\left(\sum_{u\in T}(1+\frac{w^{1}_{H_{1}}(u)}{\rho_{1}}+\frac{w^{1}_{H_{2}}(u)}{\rho_{2}}+\sum^{\infty}_{\ell=1}\sum_{0\leq i\leq \ell+1}\frac{w^{i}_{H_{1}}(u)w^{\ell+1-i}_{H_{2}}(u)}{\rho^{i}_{1}\rho^{\ell+1-i}_{2}})\right)\\
&-\left(\sum_{u\in T}\sum^{\infty}_{\ell=1}\sum_{1\leq i\leq \ell}\frac{w^{i}_{H_{1}}(u)w^{\ell+1-i}_{H_{2}}(u)}{\rho^{i}_{1}\rho^{\ell+1-i}_{2}}\right)\\
&=\sum_{u\in T}\left(1+\frac{w^{1}_{H_{1}}(u)}{\rho_{1}}+\frac{w^{1}_{H_{2}}(u)}{\rho_{2}}+(\sum^{\infty}_{\ell=1}
\frac{w^{\ell+1}_{H_{1}}(u)}{\rho^{\ell+1}_{1}}+
\frac{w^{\ell+1}_{H_{2}}(u)}{\rho^{\ell+1}_{2}})\right)\\
&=|T|+\sum^{\infty}_{\ell=1}(\frac{W^{\ell}(H_{1})}{\rho^{\ell}_{1}}+
\frac{W^{\ell}(H_{2})}{\rho^{\ell}_{2}})\\
&\geq|T|.
\end{aligned}\notag
\end{equation}
This finishes the proof of Claim 1. \hfill$\Box$

\medskip

Clearly, there are at most $2^{\frac{1}{2}|T|^{2}}$ graphs with vertex set $T$. Thus there are at most $\frac{1}{2}(2^{\frac{1}{2}|T|^{2}})^{2}=2^{|T|^{2}-1}$ pairs of graphs $Q_{1},Q_{2}$ with vertex set $T$ such that $Q_{1}\succ Q_{2}$. Let $\ell_{Q_{1},Q_{2}}$ be the integer such that $W^{\ell_{Q_{1},Q_{2}}}(Q_{1})>W^{\ell_{Q_{1},Q_{2}}}(Q_{2})$ and $W^{i}(Q_{1})=W^{i}(Q_{2})$ for any $1\leq i\leq \ell_{Q_{1},Q_{2}}-1$. Let $C=\max_{Q_{1}\succ Q_{2},V(Q_{1})=V(Q_{2})=T}\ell_{Q_{1},Q_{2}}$. Note that $C$ is only related to $|T|$. When $\rho(G)$ is sufficiently large, we can assume that $\rho(G)>|T|-1+(|T|-1)^{C+1}|T|$ in the following discussion.

$(i)$ Let $H_{1}\equiv H_{2}$. Then $W^{i}(H_{1})=W^{i}(H_{2})$ for any $i\geq1$. Without loss of generality, assume $\rho_{1}\geq\rho_{2}$. Recall $g(\rho_{1},\rho_{2})>0$. Then by Claim 1 we have
$$(\rho_{1}-\rho_{2})g(\rho_{1},\rho_{2})=\sum^{\infty}_{\ell=0}(\frac{W^{\ell+1}(H_{1})}{\rho_{1}^{\ell}}-\frac{W^{\ell+1}(H_{2})}{\rho_{2}^{\ell}})\leq0.$$
It follows that $\rho_{1}\leq\rho_{2}$. Consequently, $\rho_{1}=\rho_{2}$.

$(ii)$ Let $H_{1}\succ H_{2}$. There is an integer $1\leq k\leq C$ (defined as above), such that $W^{k}(H_{1})>W^{k}(H_{2})$ and $W^{i}(H_{1})=W^{i}(H_{2})$ for any $1\leq i\leq k-1$. We shall prove $\rho_{1}>\rho_{2}$. Suppose $\rho_{1}\leq\rho_{2}$ by contradiction. Then by Claim 1 we have
\begin{equation}
\begin{aligned}
(\rho_{1}-\rho_{2})g(\rho_{1},\rho_{2})
&=\sum^{\infty}_{\ell=0}(\frac{W^{\ell+1}(H_{1})}{\rho_{1}^{\ell}}-\frac{W^{\ell+1}(H_{2})}{\rho_{2}^{\ell}})\\
&\geq\sum^{\infty}_{\ell=0}(\frac{W^{\ell+1}(H_{1})}{\rho_{2}^{\ell}}-\frac{W^{\ell+1}(H_{2})}{\rho_{2}^{\ell}})\\
&=\sum^{\infty}_{\ell=k-1}(\frac{W^{\ell+1}(H_{1})}{\rho_{2}^{\ell}}-\frac{W^{\ell+1}(H_{2})}{\rho_{2}^{\ell}})\\
&\geq\frac{1}{\rho^{k-1}_{2}}-\sum^{\infty}_{\ell=k}\frac{W^{\ell+1}(H_{2})}{\rho_{2}^{\ell}}\\
&\geq\frac{1}{\rho^{k-1}_{2}}-\sum^{\infty}_{\ell=k}\frac{(|T|-1)^{\ell+1}|T|}{\rho_{2}^{\ell}}\\
&=\frac{1}{\rho^{k-1}_{2}(\rho_{2}-|T|+1)}(\rho_{2}-|T|+1-(|T|-1)^{k+1}|T|)\\
&\geq\frac{1}{\rho^{k-1}_{2}(\rho_{2}-|T|+1)}(\rho_{2}-|T|+1-(|T|-1)^{C+1}|T|)>0.
\end{aligned}\notag
\end{equation}
It follows that $\rho_{1}>\rho_{2}$, which contradicts the assumption $\rho_{1}\leq\rho_{2}$. Hence we have $\rho_{1}>\rho_{2}$.

$(iii)$ The proof is very similar to $(ii)$. This completes the proof. \hfill$\Box$

\section{Applications}

To characterize the graphs in ${\rm SPEX}(n,\mathcal{H})$ in Theorem \ref{Wang Kang Xue}, we first introduce some notations. For $r\geq2$ and sufficiently large $n$, let $T(n,r)$ be the Tur\'{a}n graph of order $n$ with $r$-parts $V_{1},V_{2},...,V_{r}$, where $n_{i}=|V_{i}|$ for $1\leq i \leq r$ and $\lfloor\frac{n}{r}\rfloor=n_{1}\leq n_{2}\leq\cdots\leq n_{r}=\lceil\frac{n}{r}\rceil$. For a fixed integer $t\geq1$, let $\mathcal{T}_{n,r,t}$ be the set of all the graphs obtained from $T(n,r)$ by embedding $t$ edges into the $r$-parts $V_{1},V_{2},...,V_{r}$. Given a graph $G\in\mathcal{T}_{n,r,t}$, let $H^{G}_{i}$ (probably empty) be the subgraph of $G$ induced by the edges embedded into $V_{i}$ for each $1\leq i \leq r$. Clearly, $H^{G}_{1},H^{G}_{2},...,H^{G}_{r}$ are vertex-disjoint. Let $H^{G}=\cup_{1\leq i\leq r}H^{G}_{i}$. Then $|E(H^{G})|=\sum_{1\leq i\leq r}|E(H^{G}_{i})|=t$.

\medskip

For a subset $\mathcal{G}$ of  $\mathcal{T}_{n,r,t}$, let $\mathcal{G}_{j}=\left\{G\in\mathcal{G}~|~H^{G}=\cup_{1\leq i\leq r}H^{G}_{i}=H^{G}_{j}\right\}$ for any $1\leq j\leq r$.
Define $\mathcal{H}^{\mathcal{G}}=\left\{H^{G}~|~G\in\mathcal{G}\right\}$ and $\mathcal{H}^{\mathcal{G}_{j}}=\left\{H^{G}~|~G\in\mathcal{G}_{j}\right\}$ (if $\mathcal{G}_{j}$ is not empty) for $1\leq j\leq r$.
The subset $\mathcal{G}$ is called a {\em normal} set of $\mathcal{T}_{n,r,t}$, if it satisfies the  property: for any $1\leq i\neq j\leq r$ and any $G\in\mathcal{G}_{i}$, there is a $G'\in\mathcal{G}_{j}$ such that $H^{G'}=H^{G'}_{j}$ is a copy of $H^{G}=H^{G}_{i}$.

\begin{theorem}\label{Turan}
For $t\geq1,r\geq2$ and sufficiently large $n$, let $\mathcal{G}$ be a normal subset of $\mathcal{T}_{n,r,t}$. Let $\mathcal{G}_{i},\mathcal{H^{\mathcal{G}}}$ and $\mathcal{H}^{\mathcal{G}_{i}}$ be defined as above for any $1\leq i\leq r$. If $\left\{G\in\mathcal{G}~|~H^{G}\in {\rm EX}^{2}(\mathcal{H}^{\mathcal{G}})\right\}\subseteq\cup_{1\leq i\leq r}\mathcal{G}_{i}$, then \\
$(i)$ ${\rm SPEX}(\mathcal{G})=\cup_{i,n_{i}=\lfloor\frac{n}{r}\rfloor}{\rm SPEX}(\mathcal{G}_{i})$;\\
 $(ii)$ for any $1\leq i\leq r$ and any graph $G\in\mathcal{G}_{i}$, $G\in{\rm SPEX}(\mathcal{G}_{i})$  if and only if  $H^{G}\in {\rm EX}^{\infty}(\mathcal{H}^{\mathcal{G}_{i}})$.
\end{theorem}

\f{\bf Proof:} For any $G\in\mathcal{G}$, we have $\rho(G)\geq\rho(T(n,r))$ by Lemma \ref{subgraph}, and $\rho(G)\leq\rho(T(n,r))+\rho(H^{G})$ by Courant-Weyl inequalities (see Theorem 1.3.15 of \cite{CRS}). As we know, $\rho(H^{G})\leq\sqrt{2|E(H^{G})|}=\sqrt{2t}$. Thus $\rho(T(n,r))\leq\rho(G)\leq\rho(T(n,r))+\sqrt{2t}$. Note that $\frac{\rho(T(n,r))}{n}\rightarrow1-\frac{1}{r}$ for $n\rightarrow\infty$. Thus for any $1\leq s\leq r$, we have $\frac{n_{s}}{\rho(G)}\rightarrow\frac{\frac{1}{r}}{1-\frac{1}{r}}=\frac{1}{r-1}$ for $n\rightarrow\infty$. Let $\rho=\rho(G)$.
Then by Theorem \ref{multi-set} we have
 $$\sum_{1\leq s\leq r}\frac{1}{1+\frac{n_{s}}{\rho}+\sum^{\infty}_{i=1}\frac{W^{i}(H^{G}_{s})}{\rho^{i+1}}}=r-1.$$

Let $$f(G,x)=\sum_{1\leq s\leq r}\frac{1}{1+\frac{n_{s}}{x}+\sum^{\infty}_{i=1}\frac{W^{i}(H^{G}_{s})}{x^{i+1}}},$$
 where $x>0$ is large enough such that $\sum^{\infty}_{i=1}\frac{W^{i}(H^{G}_{s})}{x^{i+1}}$ exists for $1\leq s\leq r$. Clearly, $f(G,x)$ is strictly increasing with respective to $x$. Thus $\rho$ is the largest root of $f(G,x)=r-1$.
Note that $\frac{1}{1+x}=1+\sum_{i\geq1}(-x)^{i}$ for $0\leq x<1$.
Thus for $1\leq s\leq r$, (noting that $W^{i}(H^{G}_{s})\leq 2t\cdot t^{i}$ for any $i\geq1$) we have
\begin{equation}
\begin{aligned}
&\frac{1}{1+\frac{n_{s}}{\rho}+\sum^{\infty}_{i=1}\frac{W^{i}(H^{G}_{s})}{\rho^{i+1}}}\\
&=\frac{1}{1+\frac{n_{s}}{\rho}}\frac{1}{1+\frac{1}{1+\frac{n_{s}}{\rho}}\sum^{\infty}_{i=1}\frac{W^{i}(H^{G}_{s})}{\rho^{i+1}}}\\
&=\frac{1}{1+\frac{n_{s}}{\rho}}(1-\frac{1}{1+\frac{n_{s}}{\rho}}\frac{W^{1}(H^{G}_{s})}{\rho^{2}}
-\frac{1}{1+\frac{n_{s}}{\rho}}\frac{W^{2}(H^{G}_{s})}{\rho^{3}}+\Theta(\frac{1}{\rho^{4}})).
\end{aligned}\notag
\end{equation}
Then
 $$f(G,x)=\left(\sum_{1\leq s\leq r}(\frac{1}{1+\frac{n_{s}}{x}}-\frac{1}{(1+\frac{n_{s}}{x})^{2}}\frac{W^{1}(H^{G}_{s})}{x^{2}}
-\frac{1}{(1+\frac{n_{s}}{x})^{2}}\frac{W^{2}(H^{G}_{s})}{x^{3}})\right)+\Theta(\frac{1}{x^{4}}).$$

Without loss of generality, assume $G\in\rm SPEX(\mathcal{G})$. Let $G'\in\mathcal{G}$ such that $H^{G'}\in \rm EX^{2}(\mathcal{H}^{\mathcal{G}})$. Then $G'\in\cup_{1\leq i\leq r}\mathcal{G}_{i}$ by the condition. Since $\mathcal{G}$ is a normal subset of $\mathcal{T}_{n,r,t}$, we can assume that $G'\in\mathcal{G}_{1}$.
 Note that $\rho=\rho(G)\geq\rho(G')$. Since $f(G',x)$ is strictly increasing with respective to $x$ and $f(G',\rho(G'))=r-1$, we have $f(G',\rho)\geq r-1$. It follows that $f(G',\rho)\geq f(G,\rho)$ as $f(G,\rho)=r-1$. Thus $\rho^{3}(f(G',\rho)-f(G,\rho))\geq0$. On the other hand,
 we have
\begin{equation}
\begin{aligned}
&\rho^{3}(f(G',\rho)-f(G,\rho))\\
&=\rho^{3}\left(\sum_{1\leq s\leq r}(\frac{1}{1+\frac{n_{s}}{\rho}}-\frac{1}{(1+\frac{n_{s}}{\rho})^{2}}\frac{W^{1}(H^{G'}_{s})}{\rho^{2}}
-\frac{1}{(1+\frac{n_{s}}{\rho})^{2}}\frac{W^{2}(H^{G'}_{s})}{\rho^{3}})\right)+\Theta(\frac{1}{\rho})\\
&-\rho^{3}\left(\sum_{1\leq s\leq r}(\frac{1}{1+\frac{n_{s}}{\rho}}-\frac{1}{(1+\frac{n_{s}}{\rho})^{2}}\frac{W^{1}(H^{G}_{s})}{\rho^{2}}
-\frac{1}{(1+\frac{n_{s}}{\rho})^{2}}\frac{W^{2}(H^{G}_{s})}{\rho^{3}})\right)-\Theta(\frac{1}{\rho})\\
&=\left(\rho\sum_{1\leq s\leq r}\frac{W^{1}(H^{G}_{s})-W^{1}(H^{G'}_{s})}{(1+\frac{n_{s}}{\rho})^{2}}\right)+\left(\sum_{1\leq s\leq r}\frac{W^{2}(H^{G}_{s})-W^{2}(H^{G'}_{s})}{(1+\frac{n_{s}}{\rho})^{2}}\right)+\Theta(\frac{1}{\rho}).
\end{aligned}\notag
\end{equation}

Since $G'\in\mathcal{G}_{1}$, we have that $H^{G'}=H^{G'}_{1}$, i.e., $H^{G'}_{i}$ is empty for $2\leq i\leq r$. Thus $W^{1}(H^{G'}_{1})=\sum_{1\leq s\leq r}W^{1}(H^{G}_{s})=2t$, implying that $\sum_{1\leq s\leq r}\frac{W^{1}(H^{G}_{s})-W^{1}(H^{G'}_{s})}{(1+\frac{n_{s}}{\rho})^{2}}\leq0$ (recall $n_{1}=\lfloor\frac{n}{r}\rfloor$).
From $\rho^{3}(f(G',\rho)-f(G,\rho))\geq0$, we obtain that
 $$\left(\sum_{1\leq s\leq r}\frac{W^{2}(H^{G}_{s})-W^{2}(H^{G'}_{s})}{(1+\frac{n_{s}}{\rho})^{2}}\right)+\Theta(\frac{1}{\rho})\geq0.$$
  Note that
$W^{2}(H^{G'}_{1})\geq\sum_{1\leq s\leq r}W^{2}(H^{G}_{s})=W^{2}(H^{G})$ as $H^{G'}\in \rm EX^{2}(\mathcal{H}^{\mathcal{G}})$.
If $W^{2}(H^{G'}_{1})>\sum_{1\leq s\leq r}W^{2}(H^{G}_{s})$,
then for $n\rightarrow\infty$,
 $$\sum_{1\leq s\leq r}\frac{W^{2}(H^{G}_{s})-W^{2}(H^{G'}_{s})}{(1+\frac{n_{s}}{\rho})^{2}}\leq-
\frac{1}{(1+\frac{n_{1}}{\rho})^{2}}\rightarrow-\frac{1}{(1+\frac{1}{r-1})^{2}}.$$
This contradicts that $\left(\sum_{1\leq s\leq r}\frac{W^{2}(H^{G}_{s})-W^{2}(H^{G'}_{s})}{(1+\frac{n_{s}}{\rho})^{2}}\right)+\Theta(\frac{1}{\rho})\geq0$.
   Hence $$W^{2}(H^{G'}_{1})=\sum_{1\leq s\leq r}W^{2}(H^{G}_{s})=W^{2}(H^{G}).$$
Recall that $H^{G'}\in \rm EX^{2}(\mathcal{H}^{\mathcal{G}})$. Thus $H^{G}\in \rm EX^{2}(\mathcal{H}^{\mathcal{G}})$, implying that $G\in\cup_{1\leq i\leq r}\mathcal{G}_{i}$ by the condition of the theorem.

Suppose that $G\in\mathcal{G}_{j}$ for some $1\leq j\leq r$ and $n_{j}=\lceil\frac{n}{r}\rceil>\lfloor\frac{n}{r}\rfloor$.
Then for $n\rightarrow\infty$, $$\rho\sum_{1\leq s\leq r}\frac{W^{1}(H^{G}_{s})-W^{1}(H^{G'}_{s})}{(1+\frac{n_{s}}{\rho})^{2}}=t\rho(\frac{1}{(1+\frac{n_{j}}{\rho})^{2}}-\frac{1}{(1+\frac{n_{1}}{\rho})^{2}})
=\frac{-t(2+\frac{n_{j}}{\rho}+\frac{n_{1}}{\rho})}{(1+\frac{n_{j}}{\rho})^{2}(1+\frac{n_{1}}{\rho})^{2}}$$
$$\rightarrow\frac{-t(2+\frac{2}{r-1})}{(1-\frac{1}{r-1})^{4}}.$$
Clearly, $\sum_{1\leq s\leq r}\frac{W^{2}(H^{G}_{s})-W^{2}(H^{G'}_{s})}{(1+\frac{n_{s}}{\rho})^{2}}\leq0$.
Thus for large $n$,
$$\left(\rho\sum_{1\leq s\leq r}\frac{W^{1}(H^{G}_{s})-W^{1}(H^{G'}_{s})}{(1+\frac{n_{s}}{\rho})^{2}}\right)+\left(\sum_{1\leq s\leq r}\frac{W^{2}(H^{G}_{s})-W^{2}(H^{G'}_{s})}{(1+\frac{n_{s}}{\rho})^{2}}\right)+\Theta(\frac{1}{\rho})<0.$$
Recall that $\rho^{3}(f(G',\rho)-f(G,\rho))\geq0$, i.e.,
$$\left(\rho\sum_{1\leq s\leq r}\frac{W^{1}(H^{G}_{s})-W^{1}(H^{G'}_{s})}{(1+\frac{n_{s}}{\rho})^{2}}\right)+\left(\sum_{1\leq s\leq r}\frac{W^{2}(H^{G}_{s})-W^{2}(H^{G'}_{s})}{(1+\frac{n_{s}}{\rho})^{2}}\right)+\Theta(\frac{1}{\rho})\geq0.$$
This is a contradiction. Hence $G\in\mathcal{G}_{j}$ for some $1\leq j\leq r$ and $n_{j}=\lfloor\frac{n}{r}\rfloor$.
 Thus we obtain
  $${\rm SPEX}(\mathcal{G})=\cup_{i,n_{i}=\lfloor\frac{n}{r}\rfloor}{\rm SPEX}(\mathcal{G}_{i}).$$

   For a graph $G''\in\mathcal{G}_{i}$ with $1\leq i\leq r$, note that $H^{G''}=H^{G''}_{i}$. Since $H^{G''}_{i}$ has $t$ edges, we can assume that the edges of $H^{G''}_{i}$ are embedded into a fixed subset  of size $2t$ in $V_{i}$ (this fixed subset is viewed as $T$ in Theorem \ref{one-set}).
By Theorem \ref{one-set}, we see that $G''\in{\rm SPEX}(\mathcal{G}_{i})$ if and only if  $H^{G''}=H^{G''}_{i}\in {\rm EX}^{\infty}(\mathcal{H}^{\mathcal{G}_{i}})$. This completes the proof . \hfill$\Box$

\medskip

Recall that $W^{2}(G)=\sum_{v\in V(G)}d^{2}_{G}(u)$ for any graph $G$. The following result is taken from Lemma 2.4 of \cite{IS} (it is trivial for $m=1$).

\begin{lem}{\rm (Ismailescu and Stefanica \cite{IS})}\label{2-degree}
Let $\mathcal{M}_{n,m}$ be the set of all the graphs of order $n$ with $m$ edges, where $m\geq1$ and $n\geq m+2$. Then ${\rm EX}^{2}(\mathcal{M}_{n,m})=\left\{K_{1,3}\cup \overline{K_{n-4}},K_{3}\cup \overline{K_{n-3}}\right\}$ for $m=3$, and ${\rm EX}^{2}(\mathcal{M}_{n,m})=\left\{K_{1,m}\cup \overline{K_{n-1-m}}\right\}$ otherwise.
\end{lem}

It is easy to see that $W^{3}(K_{3})=24>18=W^{3}(K_{1,3})$. Since a graph of order less than $n$ can be viewed as a graph of order $n$ by adding some isolated vertices, the following result can be deduced from Lemma \ref{2-degree} directly (by letting $n=2m+2$).

\begin{cor}\label{2-infi-walk}
Let $\mathcal{M}_{m}$ be the set of graphs with $m$ edges and without isolated vertices, where $m\geq1$. Then we have the following conclusions.\\
$(i)$ ${\rm EX}^{2}(\mathcal{M}_{m})=\left\{K_{1,3},K_{3}\right\}$ for $m=3$, and ${\rm EX}^{2}(\mathcal{M}_{m})=\left\{K_{1,m}\right\}$ otherwise\\
$(ii)$  ${\rm EX}^{\infty}(\mathcal{M}_{m})={\rm EX}^{3}(\mathcal{M}_{m})=\left\{K_{3}\right\}$ for $m=3$, and ${\rm EX}^{\infty}(\mathcal{M}_{m})={\rm EX}^{3}(\mathcal{M}_{m})=\left\{K_{1,m}\right\}$ otherwise.
\end{cor}

For a graph $G$, the {\em maximum average degree} of $G$ is the maximum value of $\frac{2|E(H)|}{|H|}$, where $H$ is an induced subgraph of $G$ with $|H|\geq1$. For $k\geq1$ and $n\geq2k+1$, Let $\mathcal{N}_{n,k}$ be the set of graphs of order $n$ with maximum average degree at most $2k$. Zhang \cite{Z} proved that any graph $G\in{\rm SPEX}(\mathcal{N}_{n,k})$ is obtained from $K_{k}\vee \overline{K_{n-k}}$ by embedding $\frac{k(k+1)}{2}$ edges into the part $V(\overline{K_{n-k}})$. When $k$ is given and $n$ is sufficiently large, from Theorem \ref{one-set} together with Corollary \ref{2-infi-walk}, we see that such $G$ is unique, of which the subgraph induced by the $\frac{k(k+1)}{2}$ edges is $K_{3}$ for $k=2$, and is $K_{1,\frac{k(k+1)}{2}}$ otherwise. (Note that, we can assume that the $\frac{k(k+1)}{2}$ edges are embedded into a fixed subset of size $k(k+1)$ which is viewed as $T$ in Theorem \ref{one-set}.)

\medskip

A graph $F$ is called color-critical if there is an edge $e$ whose deletion induces a subgraph with less chromatic number than $F$.  As in \cite{FTZ}, for $k\geq2$ color-critical graphs with chromatic number $r+1\geq3$: $F_{1},F_{2},...,F_{k}$, let $\mathcal{G}(F_{1},F_{2},...,F_{k})$ denote the family of graphs which consists of $k$ edge-disjoint copies of $F_{1},F_{2},...,F_{k}$. Fang, Tait and Zhai \cite{FTZ} proved that for sufficiently large $n$, ${\rm EX}(n,\mathcal{G}(F_{1},F_{2},...,F_{k}))$ is precisely $\mathcal{T}_{n,r,k-1}$ defined as above. From Theorem \ref{Wang Kang Xue} we see ${\rm SPEX}(n,\mathcal{G}(F_{1},F_{2},...,F_{k}))\subseteq {\rm EX}(n,\mathcal{G}(F_{1},F_{2},...,F_{k}))=\mathcal{T}_{n,r,k-1}$. The authors in
  \cite{FTZ} also devoted themselves to characterize the graphs with the maximum spectral radius in $\mathcal{T}_{n,r,k-1}$ (this result is given by Lin, Zhai and Zhao \cite{LZZ} for the case $F_{1}=F_{2}=\cdots=F_{k}=K_{3}$). In detail, they \cite{FTZ} proved the following Corollary \ref{T(n,r,k-1)}, which is directly deduced from Theorem \ref{Turan} together with Corollary \ref{2-infi-walk}. (The condition in Theorem \ref{Turan} holds for $\mathcal{G}=\mathcal{T}_{n,r,k-1}$. Clearly, $\mathcal{T}_{n,r,k-1}$ itself is a normal subset of $\mathcal{T}_{n,r,k-1}$. Moreover, the condition $$\left\{G\in\mathcal{T}_{n,r,k-1}~|~H^{G}\in {\rm EX}^{2}(\mathcal{H}^{\mathcal{T}_{n,r,k-1}})\right\}\subseteq\cup_{1\leq i\leq r}\mathcal{G}_{i}$$ holds, since the one or two graphs in ${\rm EX}^{2}(\mathcal{H}^{\mathcal{T}_{n,r,k-1}})$ are connected  by Corollary \ref{2-infi-walk}.)

\begin{cor}{\rm (\cite{FTZ,LZZ})}\label{T(n,r,k-1)}
For $k\geq2,r\geq2$ and sufficiently large $n$, let $\mathcal{T}_{n,r,k-1}$ be defined as above. Then each graph in ${\rm SPEX}(\mathcal{T}_{n,r,k-1})$ is obtained from $T_{n,r}$ by embedding $k-1$ edges into its one part of size $\lfloor\frac{n}{r}\rfloor$, of which the subgraph induced by the $k-1$ edges is $K_{3}$ for $k=4$ and is $K_{1,k-1}$ otherwise.
\end{cor}

\medskip

\f{\bf Note.} The main purpose of this paper is to supply some tools (in terms of walks) for the study of spectral radius of graphs. Thus, it will be updated in a certain time.

\medskip

\medskip

\f{\bf Acknowledgements} In the original proof of Claim 1 of Theorem \ref{one-set}, some perturbed terms are lost. Some researchers pointed out this for me and suggested to submit an updated version. I am grateful for their care and many helpful suggestions.

\medskip

\f{\bf Declaration of competing interest}

\medskip

There is no conflict of interest.

\medskip

\f{\bf Data availability statement}

\medskip

No data was used for the research described in the article.


\medskip

\begin{thebibliography}{99}









\bibitem{BDT}
 J. Byrne, D.N. Desai and M. Tait, A general theorem in spectral extremal graph theory,
arXiv:2401.07266v1.

\bibitem{CRS}
D. Cvetkovi\'{c}, P. Rowlinson and S. Simi\'{c}, An Introduction to the Theory of Graph Spectra, Cambridge University Press, Cambridge, 2010.

\bibitem{CFTZ}
S. Cioab\u{a}, L. Feng, M. Tait and X. Zhang, The maximum spectral radius of graphs
without friendship subgraphs, Electron. J. Combin. 27 (2020)  \#4.22.

\bibitem{CDT1}
S. Cioab\u{a}, D. Desai and M. Tait, The spectral even cycle problem, arXiv:2205.00990v1.


\bibitem{CDT2}
S. Cioab\u{a}, D. Desai and M. Tait, The spectral radius of graphs with no odd wheels,
European J. Combin. 99 (2022) 103420.

\bibitem{CDT3}
S. Cioab\u{a}, D. Desai and M. Tait,
A spectral Erd\H{o}s-S\'{o}s theorem, SIAM J. Discrete
Math. 37 (2023) 2228-2239.

\bibitem{DKLNTW}
 D. Desai, L. Kang, Y. Li, Z. Ni, M. Tait and J. Wang, Spectral extremal graphs
for intersecting cliques, Linear Algebra Appl. 644 (2022) 234-258.

\bibitem{FLSZ}
L. Fang, H. Lin, J. Shu and Z. Zhang, Spectral extremal results on trees, arXiv:2401.05786v2.
\bibitem{FTZ}
L. Fang, M. Tait and M. Zhai, Tur\'{a}n numbers for non-bipartite graphs and applications to spectral extremal problems, arXiv:2404.09069v1.

\bibitem{IS}
D. Ismailescu and D. Stefanica, Minimizer graphs for a class of extremal problems, J. Graph Theory 39 (2002) 230-240.

\bibitem{LL}
 X. Lei and S. Li, Spectral extremal problem on disjoint color-critical graphs, Electron. J. Combin. 31 (2024) \# 1.25.

\bibitem{LP1}
Y. Li, and Y. Peng, Refinement on spectral Tur\'{a}n's theorem, SIAM J. Discrete Math.
37 (2023) 2462-2485.
\bibitem{LP2}
Y. Li, and Y. Peng, The spectral radius of graphs with no intersecting odd cycles,
Discrete Math. 345 (2022) 112907.

\bibitem{LZZ}
 H. Lin, M. Zhai and Y. Zhao, Spectral radius, edge-disjoint cycles and cycles of
the same length, Electron. J. Combin. 29 (2022) \# 2.1.

\bibitem{NWK}
 Z. Ni, J. Wang and L. Kang, Spectral extremal graphs for disjoint cliques, Electron.
J. Combin. 30 (2023)  \#1.20.

\bibitem{N1}
V. Nikiforov, The spectral radius of graphs without paths and cycles of specified length, Linear Algebra Appl. 432 (2010) 2243-2256.

\bibitem{N2}
 V. Nikiforov, Bounds on graph eigenvalues II, Linear Algebra Appl. 427 (2007) 183-189.


\bibitem{WKX}
J. Wang, L. Kang and Y. Xue, On a conjecture of spectral extremal problems, J.
Combin. Theory Ser. B 159 (2023) 20-41.

\bibitem{ZL}
 M. Zhai and H. Lin, Spectral extrema of graphs: Forbidden hexagon, Discrete Math.
343 (2020) 112028.

\bibitem{Z}
W. Zhang, The spectral radius, maximum average degree and cycles
of consecutive lengths of graphs, Graphs Combin. 40 (2024) 32.



\end{thebibliography}
\end{document}